\documentclass{svjour3}    
\usepackage[utf8]{inputenc}
\usepackage{amsmath}
\usepackage{graphicx}
\usepackage{algorithm,algorithmicx,algpseudocode}
\usepackage{url}

\title{Accurate algorithms for Bessel matrices}

\author{Jorge Delgado, H\'ector Orera, J. M. Pe\~na}

\institute{J. Delgado \at
Departamento de Matem\'atica Aplicada \\ 
Escuela Universitaria Politécnica de Teruel \\
Universidad de Zaragoza \\ 
44003 Teruel, Spain. \\
              Tel.: +34-978618174\\
              \email{jorgedel@unizar.es}           %  \\
           \and
           H. Orera \at
             Departamento de Matem\'atica Aplicada \\
Facultad de Ciencias \\ 
             Universidad de Zaragoza \\ 
             50009 Zaragoza, Spain
                        \and
           J.~M. Pe\~na \at
             Departamento de Matem\'atica Aplicada \\
Facultad de Ciencias \\ 
             Universidad de Zaragoza \\ 
             50009 Zaragoza, Spain
}

\date{Received: date / Accepted: date}

\begin{document}

\maketitle

\begin{abstract}
In this paper, we prove that any collocation matrix of Bessel polynomials at positive points is strictly totally positive, that is, all its minors are positive. Moreover, an accurate method to construct the bidiagonal factorization of these matrices is obtained and used to compute with high relative accuracy the eigenvalues, singular values and inverses.  Similar results for the collocation matrices for the reverse Bessel polynomials are also obtained. Numerical examples illustrating the theoretical results are included.

\end{abstract}

\section{Introduction}

Bessel polynomials are ubiquitous and occur in many fields such as partial differential equations, number theory, algebra and statistics (see \cite{G}). They  form an orthogonal sequence of polynomials and are related to the modified Bessel function of the second kind (see pp. 7 and 34 of \cite{G}). They are closely related to the reverse Bessel polynomials, which have many applications in Electrical Engineering. In particular, they play a key role in network analysis of electrical circuits 
(see page 145 of \cite{G} and references therein).
In Combinatorics, the coefficients of the reverse Bessel polynomials are also known as signless Bessel numbers of the first kind. The Bessel numbers have been studied from a combinatorial perspective and are closely related to the Stirling numbers \cite{H,Y}. 
In \cite{Kr} it was shown that Bessel polynomials occur naturally in the theory of traveling spherical waves. 
 Bessel polynomials are also very important for some problems of static potentials, signal processing and electronics. For example, the Bessel polynomials are used in Frequency Modulation (FM) synthesis and in the Bessel filter. In the case of FM synthesis, the polynomials are used to compute the bandwidth of a modulated in frequency signal.
The zeros of Bessel polynomials and generalized Bessel polynomials also play a crucial role in applications in Electrical Engineering. On the accurate 
computations of the zeros of generalized Bessel polynomials see \cite{PASQ}. 

This paper deals with the accurate computation when using collocation matrices of Bessel polynomials and reverse Bessel polynomials. It is shown that these  matrices provide new  structured classes  for which algebraic computations (such as the computation of the inverse, of all the eigenvalues and singular values, or the solutions of some linear systems) can be performed with high relative accuracy (HRA). Moreover, a crucial result for this purpose has been the total positivity of the considered matrices. Let us recall that a matrix is {\it totally positive} ({\it strictly totally positive}) if all its minors are nonnegative (positive) and will be denoted TP (STP). These matrices have also been called
in the literature totally nonnegative (totally positive). Many applications of these matrices can be seen in \cite{A,FJ,p}.
For some subclasses of TP matrices a bidiagonal factorization with HRA has been obtained (cf. \cite{Koev.Dem,Pascal,JS,qBer,BV,SBV,BV2}).
In \cite{cmp} it was proved that the first positive zero of a Bessel function of the first kind is the half of the critical length of a cycloidal space, relating Bessel functions with the total positivity theory and computer-aided geometric design. We prove here a new surprising connection of  total positivity with Bessel functions through the collocation matrices of Bessel polynomials.

The paper is organized as follows. In Section 2, we present some auxiliary results and basic notations related to the bidiagonal factorization of totally positive matrices as well as with high relative accuracy. In Section 3 we introduce the Bessel polynomials and we first prove that the matrix of change of basis between the Bessel polynomials and the monomials is TP. We also define the Bessel matrices and prove that they are STP. Finally, we   provide the construction with high relative accuracy of the bidiagonal factorization of Bessel matrices, which in turn can be used to apply the algorithms presented by P. Koev in \cite{KoevSoft} for the algebraic computations mentioned above with 
high relative accuracy. A similar task for the collocation matrices of the reverse Bessel polynomials is performed in Section 4.    Finally, Section 5 includes illustrative numerical examples confirming the theoretical results for the computation of eigenvalues, singular values, inverses, and the solution of linear systems with the matrices considered in this paper.

\section{Auxiliary results}\label{sec.aux}

 Neville elimination (NE) is an alternative procedure to Gaussian elimination.
 NE makes zeros in a column of a matrix by adding to 
 each row an appropriate multiple of the previous one (see \cite{TPandNE}).
 Given a nonsingular  
matrix
$A=(a_{ij})_{1\leq i,j\leq n}$, the NE procedure consists of $n-1$  steps, leading to a
sequence of matrices as follows:
\begin{equation} \label{seq}
A=A^{(1)}\rightarrow \widetilde{A}^{(1)}\rightarrow A^{(2)}\rightarrow \widetilde{A}^{(2)}
        \rightarrow \cdots \rightarrow A^{(n)}= \widetilde{A}^{(n)}=U,
\end{equation}
with $U$  an upper triangular matrix.

On the one hand, $\widetilde{A}^{(t)}$ is obtained from the matrix $A^{(t)}$ by
moving to the bottom the rows with a zero entry in column $t$ below the main
diagonal, if necessary. The matrix $A^{(t+1)}$, $t=1,\ldots,n-1$, 
is obtained from $\widetilde{A}^{(t)}$ by computing
\begin{equation} \label{elemen_col}
    a_{ij}^{(t+1)}=
    \left\{%
\begin{array}{ll}
%    a_{ij}^{(t)}, & \hbox{if } 1\leq i\leq t, \\
    \widetilde{a}_{ij}^{(t)}-\displaystyle{\frac{\widetilde{a}_{it}^{(t)}}{\widetilde{a}_{i-1,t}^{(t)}}}
    \widetilde{a}_{i-1,j}^{(t)}, &
         \hbox{if } t\leq j< i\leq n\hbox{ and } \widetilde{a}_{i-1,t}^{(t)}\neq 0, \\
%    a_{ij}^{(t)}, & \hbox{if } t\leq j< i\leq n\hbox{ and } a_{i-1,t}^{(t)}=0, \\
    \widetilde{a}_{ij}^{(t)}, & \hbox{otherwise,} 
\end{array}%
\right.
\end{equation}
for all $t \in \{1,\dots,n-1\}$.

The entry
\begin{equation} \label{pivote}
    p_{ij}:=\widetilde{a}_{ij}^{(j)},\quad 1\leq j\leq i\leq n,
\end{equation}
is  the $(i,j)$ {\it pivot} of the NE of $A$, and the pivots $p_{ii}$ are called 
{\it diagonal pivots}. The number
\begin{equation}\label{multiplicador}
	m_{ij}=\left\{
	\begin{array}{ll}
		\frac{\widetilde{a}_{ij}^{(j)}}{\widetilde{a}_{i-1,j}^{(j)}}=\frac{p_{ij}}{p_{i-1,j}}, & 
		\text{if } \widetilde{a}_{i-1,j}^{(j)}\ne 0, \\
		0, & \text{if } \widetilde{a}_{i-1,j}^{(j)}=0,
	\end{array}
	\right.
\end{equation}
is called the $(i,j)$ {\it multiplier} of NE 
of $A$, where $1\leq j<i\leq n$.

NE is a very useful method when applied to TP matrices. If $A$ is a nonsingular TP matrix, then no rows exchanges are needed
when applying NE and so, in this 
case, $A^{(t)}=\widetilde{A}^{(t)}$ for all $t$. In fact, in Theorem 5.4 of \cite{TPandNE} the following characterization of nonsingular TP matrices was provided.

\begin{theorem}\label{thm:TP.car}
Let $A$ be a nonsingular matrix. Then $A$ is TP if and only if there are no row exchanges in the NE of $A$ and $U^T$ and the pivots of both NE are  nonnegative.
\end{theorem}

In \cite{fact} it was seen that nonsingular TP matrices can be expressed as a unique
bidiagonal decomposition. 

\begin{theorem}\label{BD}
 (cf. Theorem 4.2 of  \cite{fact}).
 Let $A$ be a nonsingular $n\times n$ TP matrix. 
 Then $A$ admits a decomposition of the form
 \begin{equation}\label{f}  
  A=F_{n-1}\cdots F_1 DG_1\cdots G_{n-1},
 \end{equation}
 where $F_i$ and $G_i$, $i\in\{1,\ldots,n-1\}$, are the lower and upper 
 triangular nonnegative bidiagonal matrices given by
 \begin{equation}\label{f.det}
  F_i=\left(
    \begin{smallmatrix}
     1 &        &                   &            &  &  &  \\
     0 & 1      &                    &            &  &  &  \\
       & \ddots & \ddots            &            &  &  &  \\
       &        & 0      & 1                      &  &  &  \\
       &        &                    & m_{i+1,1} & 1 &  & \\
       &        &                    &            & \ddots & \ddots & \\
       &        &                    &            &  & m_{n,n-i} & 1
    \end{smallmatrix}
  \right),\;
 G_i^T=\left(
    \begin{smallmatrix}             
     1 &        &                    &            &  &  &  \\
     0 & 1      &                    &            &  &  &  \\
       & \ddots & \ddots             &            &  &  &  \\
       &        & 0      & 1                      &  &  &  \\
       &        &                    & \widetilde{m}_{i+1,1} & 1 &  & \\
       &        &                    &            & \ddots & \ddots & \\
       &        &                    &            &  & \widetilde{m}_{n,n-i} & 1
    \end{smallmatrix}
  \right),
 \end{equation}
 and $D$ a diagonal matrix $diag(p_{11},\ldots,p_{nn})$ with positive diagonal entries. If, in addition, the entries  $m_{ij}$, $\widetilde{m}_{ij}$ satisfy
 $$m_{ij}=0\Rightarrow m_{hj}=0 \quad \forall \, h>i$$
 and
 $$\widetilde{m}_{ij}=0\Rightarrow m_{ik}=0 \quad \forall \, k>j,$$
 then the decomposition (\ref{f}) is unique.
 \end{theorem}

 By Theorems 4.1 and 4.2 of  \cite{fact} we also know  that, for $1\leq j<i\leq n$,
 $m_{ij}$ and $p_{ii}$ in the bidiagonal decomposition given by \eqref{f} with \eqref{f.det} 
 are the multipliers and the diagonal pivots
 when applying the NE to $A$ and, using the arguments of p. 116 of \cite{fact}, $\widetilde{m}_{ij}$
 are the multipliers when applying the NE to $A^T$.

 In \cite{K2} it was devised a concise matrix notation $\mathcal{BD}(A)$ 
 for the bidiagonal decomposition \eqref{f} and \eqref{f.det} given by
 \begin{equation}\label{eq:BD.A}
 	(\mathcal{BD}(A))_{ij}=\left\{
 		\begin{array}{ll}
 			m_{ij}, & \text{if } i>j, \\
 			\widetilde{m}_{ji}, & \text{if } i<j, \\
 			p_{ii}, & \text{if } i=j.
 		\end{array}
	\right.
 \end{equation}
 
 \begin{remark}\label{rem:bd}
  If $A$ is a TP matrix, then $A^T$ is also TP. 
  Transposing formula \eqref{f} of Theorem \ref{BD} we obtain 
  the unique bidiagonal decomposition of $A^T$:
  \[
  	A^T=G_{n-1}^T\cdots G_1^T D F_1^T\cdots F_{n-1}^T,
  \]
  where $F_i$ and $G_i$, $i\in\{1,\ldots,n-1\}$, are the lower and upper 
 triangular nonnegative bidiagonal matrices given in formula \eqref{f.det}. It
 can also be checked that 
 \[
 	\mathcal{BD}(A^T)=\mathcal{BD}(A)^T.
 \]
 \end{remark}

An algorithm can be performed with high relative accuracy if
all the included subtractions are of initial data,
that is, if it only includes products, 
divisions, sums of numbers of the same sign and subtractions of the initial data 
(cf. \cite{Koev.Dem,K2}). In \cite{K2}, given a nonsingular TP matrix $A$ whose $\mathcal{BD}(A)$ is known with HRA, Koev presented algorithms 
for computing to HRA the eigenvalues of the matrix $A$, the singular values of the
 matrix $A$, the inverse of the 
matrix $A$ and the solution of linear systems of equations $Ax=b$, where $b$ has a pattern of alternating signs.

\section{Bessel polynomials and matrices}

Let us consider the Bessel polynomials defined by
\begin{equation}\label{pol}
B_n(x)=\sum_{k=0}^n \frac{(n+k)!}{2^k(n-k)!k!}x^k, \quad  n=0,1,2\ldots,
\end{equation}

Given a real positive integer $n$, let us define the corresponding \textit{semifactorial} by
\begin{equation*}
n!!=\prod^{[n/2]-1}_{k=0}(n-2k).
\end{equation*}

Let $A=(a_{ij})_{1 \leq i,j\leq n}$ be the lower triangular matrix such that
\begin{equation}\label{eq:cambioM}
	(B_0(x),B_1(x),\ldots,B_{n-1}(x))^T=A(1,x,\ldots,x^{n-1})^T,
\end{equation}
that is, the lower triangular matrix $A$ is defined by 
\begin{equation}\label{eq:def.matC}
	a_{ij}:=\left\{
		\begin{array}{ll}
			\frac{(i+j-2)!}{2^{j-1}(i-j)!(j-1)!}=\frac{(2j-2)!}{2^{j-1}(j-1)!}{i+j-2 \choose i-j}, & \text{if } i\ge j, \\
			0, & \text{if } i<j.
		\end{array}
	\right.
\end{equation}

Theorem \ref{thm.bessel} proves the total positivity of $A$, and provides $\mathcal{BD}(A)$. In addition, its proof gives the explicit form of all the entries of the matrices $A^{(k)}$ of \eqref{seq} computed through the NE of $A$.

\begin{theorem}\label{thm.bessel}
Let $A=(a_{ij})_{1 \leq i,j\leq n}$ be the lower triangular matrix in \eqref{eq:cambioM} defined by \eqref{eq:def.matC}. Then we have that 
\begin{itemize}
	\item[(i)] the pivots of the NE of $A$ are given by
	\begin{equation}\label{eq:piv.A}
	p_{ij}=
	\begin{array}{ll}
		\frac{1}{2^{j-1}}\frac{(i-1)!}{(i-j)!}\prod_{r=1}^{j-1}\frac{(2i-r-1)}{(i-j+r)}, & 
	1\leq j\leq i \leq n, 
	\end{array}
\end{equation}
and the multipliers by
\begin{equation}\label{eq:mult.A}
m_{ij}=
	\begin{array}{ll}
		\frac{(2i-2)(2i-3)}{(2i-j-1)(2i-j-2)}, & 
		1\leq j < i \leq n, 
	\end{array}
\end{equation}
	\item[(ii)] $A$ is a nonsingular TP matrix
	\item[(iii)] and the bidiagonal factorization of $A$ is given by
\begin{equation}\label{BD.bessel}
\mathcal{BD}(A)_{ij}=\left\{
 		\begin{array}{ll}
 			\frac{(2i-2)(2i-3)}{(2i-j-1)(2i-j-2)}, & \text{if } i>j, \\
 			1, & \text{if } i=j=1, \\
 			(2i-3)!!, & \text{if } i=j>1,\\
 			0, & \text{if } i<j, 			
 		\end{array}
	\right.
\end{equation}
and can be computed to HRA.
\end{itemize}
\end{theorem}

\begin{proof}
 \begin{itemize}
 \item[(i)] 
If $A^{(k)}=(a^{(k)}_{ij})_{1 \leq i,j\leq n}$ is the matrix obtained after $k-1$ steps of the NE of $A$ (see \eqref{seq})
for $k=2,\ldots,n$, let us prove by induction on $k\in\{2,\ldots,n\}$ that
\begin{equation}\label{NE_A}
a^{(k)}_{ij}=\frac{1}{2^{j-1}}\frac{(i+j-k-1)!}{(i-j)!(j-k)!}\prod_{r=1}^{k-1}\frac{(2i-r-1)}{(i-k+r)}
\end{equation} 
for $i\ge j$.  For the case $k=2$, let us start by computing the first step of the NE of $A$:
\begin{align*}
a^{(2)}_{ij}&=a_{ij}-\frac{a_{i1}}{a_{i-1,1}}a_{i-1,j}=a_{ij}-a_{i-1,j}\\
&=\frac{1}{2^{j-1}}\left(\frac{(i+j-2)!}{(i-j)!(j-1)!}-\frac{(i+j-3)!}{(i-j-1)!(j-1)!}\right)\\
&=\frac{(i+j-3)!}{2^{j-1}(i-j-1)!(j-1)!}\left(\frac{i+j-2}{i-j}-1\right)\\
&=\frac{(i+j-3)!(2j-2)}{2^{j-1}(i-j)!(j-1)!}\\
&= \frac{(i+j-3)!}{2^{j-1}(i-j)!(j-2)!}\frac{(2i-2)}{(i-1)}.
\end{align*}

Hence, formula $\eqref{NE_A}$ holds for $k=2$. Now, let us assume that $\eqref{NE_A}$ holds for some $k\in\{2,\ldots,n-1\}$ and let us prove that  it also holds for $k+1$. 
Performing the $k$-th step of the NE we have
\[
a^{(k+1)}_{ij}=a^{(k)}_{ij}- \frac{a^{(k)}_{ik}}{a^{(k)}_{i-1,k}}a^{(k)}_{i-1,j}.
\]
Then, by using the induction hypothesis, we obtain
\begin{align*}
a^{(k+1)}_{ij}&=a^{(k)}_{ij}- \frac{a^{(k)}_{ik}}{a^{(k)}_{i-1,k}}a^{(k)}_{i-1,j}\\
&=a^{(k)}_{ij}- \frac{(2i-2)(2i-3)}{(2i-k-1)(2i-k-2)}a^{(k)}_{i-1,j}\\
&= \frac{(i+j-k-1)!}{2^{j-1}(i-j)!(j-k)!}\prod_{r=1}^{k-1}\frac{(2i-r-1)}{(i-k+r)}\\
&\phantom{=}- \frac{(2i-2)(2i-3)}{(2i-k-1)(2i-k-2)}\frac{(i+j-k-2)!}{2^{j-1}(i-j-1)!(j-k)!}\prod_{r=1}^{k-1}\frac{(2i-r-3)}{(i-k+r-1)}.\\
\end{align*}
Simplifying the previous formula, $a^{(k+1)}_{ij}$ can be written as 
\begin{equation*}
a^{(k+1)}_{ij}=\frac{(i+j-k-2)!}{2^{j-1}(i-j-1)!(j-k)!}\prod_{r=1}^{k-1}\frac{(2i-r-1)}{(i-k+r)}\left(\frac{i+j-k-1}{i-j} - \frac{i-1}{i-k} \right).
\end{equation*}
From the previous expression we can deduce that
\begin{align*}
a^{(k+1)}_{ij}&=\frac{(i+j-k-2)!}{2^{j-1}(i-j-1)!(j-k)!}\prod_{r=1}^{k-1}\frac{(2i-r-1)}{(i-k+r)} \cdot \frac{2ji-2ki-kj-j+k^2+k}{(i-j)(i-k)}\\
&=\frac{(i+j-k-2)!}{2^{j-1}(i-j-1)!(j-k)!}\prod_{r=1}^{k-1}\frac{(2i-r-1)}{(i-k+r)} \cdot \frac{(j-k)(2i-k-1)}{(i-j)(i-k)}\\
&=\frac{(i+j-k-2)!}{2^{j-1}(i-j)!(j-k-1)!}\prod_{r=1}^{k}\frac{(2i-r-1)}{(i-k+r-1)}.
\end{align*}
Therefore, $\eqref{NE_A}$ holds for $k+1$ 
and the result follows.
 
  The pivot  $p_{ij}=\tilde{a}^{(j)}_{ij}=a^{(j)}_{ij}$ is given by \eqref{NE_A} with $k=j$ and we have that, for $i>j$, $m_{ij}=\frac{p_{ij}}{p_{i-1,j}}$, obtaining formulas \eqref{eq:piv.A} and \eqref{eq:mult.A}, respectively. 
 \item[(ii)] The  lower triangular matrix $A$ is nonsingular
 	since it has nonzero diagonal entries. It can be seen in the proof of (i) that the NE of $A$
 satisfies the hypotheses of Theorem \ref{thm:TP.car}. Since $U^T$ is a diagonal matrix, the NE of $U^T$
 obviously satisfies the hypotheses of Theorem \ref{thm:TP.car}.
 	Hence, we can conclude that $A$ is a TP matrix.  
 \item[(iii)] By (i) and taking into account that $U$ is a diagonal matrix with diagonal entries $a_{ii}$ ($1\le i\le n$), it is straightforward to deduce that $\mathcal{BD}(A)$ is given by \eqref{BD.bessel}.
 The subtractions in this formula are of integers and, hence, they can be computed to HRA, in fact,
 in an exact way.
 \end{itemize}
\end{proof}

Let us introduce the collocation matrices of the Bessel polynomials.
\begin{definition}
	Given a sequence of parameters $0<t_0<t_1<\cdots<t_{n-1}$ we call the collocation
	matrix of the Bessel polynomials $(B_0,\ldots,B_{n-1})$ at that sequence,
	\[
	M=M\left(
		\begin{array}{c}
			B_0,\ldots,B_{n-1} \\
			t_0,\ldots,t_{n-1}
		\end{array}
	\right)=(B_{j-1}(t_{i-1}))_{1\leq i,j\leq n},
	\]
	a {\it Bessel matrix}.
\end{definition}

The following result proves that the Bessel matrices are STP and that some
usual algebraic problems with these matrices can be solved to HRA.
\begin{theorem}\label{thm:HRA.Bes}
	Given a sequence of parameters $0<t_0<t_1<\cdots<t_{n-1}$, the corresponding Bessel matrix 
	$M$ is an STP matrix and given the parametrization $t_i$ ($0\leq i \leq n-1$), the following
computations can be performed with HRA: all the eigenvalues, all the singular values, the inverse of 
the Bessel matrix $M$, and the solution of the linear systems $Mx=b$, where $b=(b_1,\ldots,b_n)^T$ has alternating signs.
\end{theorem}
\begin{proof}
By formula \eqref{eq:cambioM} we have that
\[
	M=VA^T,
\]
where $M$ is the Bessel matrix corresponding to the collocation matrix 
of the Bessel polynomials $(B_0,\ldots,B_{n-1})$ at $t_0,\ldots,t_{n-1}$,
$A$ is the lower triangular matrix defined by \eqref{eq:def.matC} and
$V$ is the Vandermonde matrix corresponding to the collocation matrix of
the monomial basis of degree $n-1$ at $t_0,\ldots,t_{n-1}$. Since $V$
is a Vandermonde matrix with strictly increasing positive nodes, it is STP
(see page 111 of \cite{Gant} and page 12 of \cite{FJ}). $A^T$ is nonsingular TP because $A$ 
is nonsingular TP by 
Theorem \ref{thm.bessel} (ii). Then, by Theorem 3.1 of \cite{A}, the Bessel matrix $M$ is 
STP because it is the product of an STP matrix and a nonsingular TP matrix.

In Section 5.2 of \cite{K2} Koev devised an algorithm (Algorithm 5.1 in the reference) 
that, given the bidiagonal decompositions $\mathcal{BD}(C)$ and $\mathcal{BD}(D)$
to HRA of two TP matrices $C$ and $D$, provides the bidiagonal decomposition 
$\mathcal{BD}(CD)$ to HRA of the TP product matrix $CD$. Since the bidiagonal factorization 
$\mathcal{BD}(V)$ of the Vandermonde matrix $V$ is known 
to HRA (see Section 3 of \cite{K1}) and the bidiagonal factorization $\mathcal{BD}(A^T)$ 
of the matrix $A$ can be computed to HRA by
Theorem \ref{thm.bessel} (iii) and Remark \ref{rem:bd}, by using Algorithm 5.1 of \cite{K2} the bidiagonal factorization
$\mathcal{BD}(M)=\mathcal{BD}(VA^T)$ of the Bessel matrix $M$ is obtained to HRA.

Finally, the construction of $\mathcal{BD}(M)$ with HRA guarantees that the algebraic 
computations mentioned in the statement of this theorem can be performed with HRA (see 
Section \ref{sec.aux} of this paper or Section 3 of \cite{K2}).
\end{proof}

A system of functions is STP when all its collocation matrices are STP. 
The following result is a straightforward consequence of the previous theorem.
\begin{corollary}
The system of functions formed by the Bessel polynomials of degree less than $n$,
$(B_0(x),B_1(x),\ldots,B_{n-1}(x))$, $x\in(0,+\infty)$, is an STP system.
\end{corollary}

\section{Reverse Bessel polynomials and matrices}

Reversing the order of the coefficients of $B_n (x)$ in \eqref{pol} we can define the {\it reverse Bessel polynomials}:

\begin{equation}\label{rev.pol}
B^r_n(x)=\sum_{k=0}^n \frac{(n+k)!}{2^k(n-k)!k!}x^{n-k}, \quad  n=0,1,2\ldots,
\end{equation}

Let $C=(c_{ij})_{1 \leq i,j\leq n}$ be the lower triangular matrix such that
\begin{equation}\label{eq:cambioMr}
	(B^r_0(x),B^r_1(x),\ldots,B^r_{n-1}(x))^T=C(1,x,\ldots,x^{n-1})^T,
\end{equation}
that is, the lower triangular matrix $C$ is defined by 
\begin{equation}\label{eq:def.matCr}
	c_{ij}=\left\{
		\begin{array}{ll}
			\frac{(2i-j-1)!}{2^{i-j}(j-1)!(i-j)!}, & i\ge j, \\
			0, & i<j.
		\end{array}
	\right.
\end{equation}

Theorem \ref{thm.rev.bessel} proves the total positivity of $C$, and provides $\mathcal{BD}(C)$. In addition, its proof gives the explicit form of all the entries of the matrices $C^{(k)}$ computed through the NE of $C$.

\begin{theorem}\label{thm.rev.bessel}
Let $C=(c_{ij})_{1 \leq i,j\leq n}$ be the lower triangular matrix in
\eqref{eq:cambioMr} defined by \eqref{eq:def.matCr}. Then, we have that 
\begin{itemize}
\item[(i)] the pivots of the NE of $C$ are given by
\begin{equation}\label{eq:piv.C}
	\begin{array}{rll}
	&p_{ij}=\frac{(2i-2j)!}{2^{i-j}(i-j)!} \quad 1\leq j\leq i\leq n &\text{if } j \text{ is odd}, \\
	&p_{ij}=0\quad 1\leq j<i\leq n,\quad p_{jj}=1\quad 1\leq j\leq n &\text{if } j \text{ is even},
	\end{array}
\end{equation}
and the multipliers by 
\begin{equation}\label{eq:mult.C}
	\begin{array}{rll}
	&m_{ij}=2i-1-2j \quad 1\leq j< i\leq n &\text{if } j \text{ is odd}, \\
	&m_{ij}=0,\quad 1\leq j< i\leq n &\text{if } j \text{ is even},
	\end{array}
\end{equation}
\item[(ii)] $C$ is a nonsingular TP matrix
\item[(iii)] and the bidiagonal factorization of $C$ is given by
\begin{equation}\label{BD.rev.bessel}
\mathcal{BD}(C)_{ij}=\left\{
 		\begin{array}{ll}
 			2i-2j-1, & \text{if } i>j \text{ with }j\text{ odd,}\\
 			1, & \text{if } i=j,\\
 			0, & \text{otherwise},
 		\end{array}
	\right.
\end{equation}
and can be computed to HRA.
\end{itemize}
\end{theorem}
\begin{proof}
\begin{itemize}
\item[(i)] Let us perform the first step of the NE of $C$:
\[
	c^{(2)}_{ij}=c_{ij}-\frac{c_{i1}}{c_{i-1,1}}c_{i-1,j}=c_{ij}-(2i-3)c_{i-1,j},\quad i>j\ge 1.
\]
By using \eqref{eq:def.matCr} in the previous expression and simplifying we have that
\begin{align*}
	c^{(2)}_{ij}&=\frac{(2i-j-1)!}{2^{i-j}(j-1)!(i-j)!}-\frac{(2i-3)(2i-j-3)!}{2^{i-j-1}(j-1)!(i-j-1)!} \\
		&=\frac{(2i-j-3)!}{2^{i-j-1}(j-1)!(i-j-1)!}\left[\frac{(2i-j-2)(2i-j-1)}{2(i-j)}-2i+3\right] \\
		&=\frac{(2i-j-3)!}{2^{i-j}(j-1)!(i-j)!}(j^2-3j+2)\\
		&=\frac{(2i-j-3)!}{2^{i-j}(j-1)!(i-j)!}(j-2)(j-1).
\end{align*}
From the previous formula we deduce that
\begin{equation}\label{eq:NE.C.1}
c_{ij}^{(2)}=\left\{
	\begin{array}{ll}
		\frac{(2i-j-3)!}{2^{i-j}(j-3)!(i-j)!}, & \text{for } i>j\ge 3, \\
		0, & \text{for } i>j\le 2, \\
		c_{ij}^{(1)}, & \text{in otherwise.}
	\end{array}
\right.
\end{equation}
Since $c_{i2}^{(2)}=0$ for $i=3,\ldots,n$, we have that $C^{(3)}=C^{(2)}$.  
So, taking into account this fact, formula \eqref{eq:NE.C.1} and formula \eqref{eq:def.matCr}, we can observe that $c^{(3)}_{ij}=c_{i-2,j-2}$  and so, performing two steps of the NE of $C$ gives as a result a leading principal submatrix of $C$. In particular, $C$ satisfies that $C^{(2)}[3,\ldots,n]=C^{(3)}[3,\ldots,n]=C[1,\ldots,n-2]$. Hence, and by formulas \eqref{pivote} and
\eqref{multiplicador}, we deduce formulas \eqref{eq:piv.C} and \eqref{eq:mult.C}.
\item[(ii)] The  lower triangular matrix $C$ is nonsingular
 	since it has nonzero diagonal entries. It can be seen in the proof of (i) that the NE of $C$
 satisfies the hypotheses of Theorem \ref{thm:TP.car}. Since the upper triangular matrix obtained after the process (\ref{seq}) of the NE of $C$ is a diagonal matrix, the NE of its transpose
 obviously satisfies the hypotheses of Theorem \ref{thm:TP.car}. Hence, we can conclude that $C$ is a TP matrix.
\item[(iii)] By (i) it is straightforward to deduce that $\mathcal{BD}(C)$ is given by \eqref{BD.rev.bessel}. 	
 The subtractions in this formula are of integers and, hence, they can be computed to HRA, in fact,
 in an exact way.
\end{itemize}
\end{proof}

\begin{definition}
	Given a sequence of parameters $0<t_0<t_1<\cdots<t_{n-1}$ we call the collocation
	matrix of the reverse Bessel polynomials $(B^r_0,\ldots,B^r_{n-1})$ at that sequence,
	\[
	M_r=M\left(
		\begin{array}{c}
			B^r_0,\ldots,B^r_{n-1} \\
			t_0,\ldots,t_{n-1}
		\end{array}
	\right)=(B_{j-1}^r(t_{i-1}))_{1\leq i,j\leq n}
	\]
	a {\it reverse Bessel matrix}.
\end{definition}

The following result proves that the reverse Bessel matrices are STP 
and that some usual algebraic problems with these matrices can be solved to HRA.

\begin{theorem}
	Given a sequence of parameters $0<t_0<t_1<\cdots<t_{n-1}$, the corresponding reverse Bessel matrix 
	$M_r$ is an STP matrix and given the parametrization $t_i$ ($0\leq i \leq n-1$), the following
computations can be performed with HRA: all the eigenvalues, all the singular values, the inverse of 
the reverse Bessel matrix $M_r$, and the solution of the linear systems $M_rx=b$, where $b=(b_1,\ldots,b_n)^T$ has alternating signs.
\end{theorem}
\begin{proof}
The results can be proved in an anologous way to those of Theorem \ref{thm:HRA.Bes}.
\end{proof}

The following result is a straightforward consequence of the previous theorem.
\begin{corollary}
The system of functions formed by the reverse Bessel polynomials of degree less than $n$,
$(B_0^r(x),B_1^r(x),\ldots,B_{n-1}^r(x))$, $x\in(0,+\infty)$, is an STP system.
\end{corollary}

\section{Numerical experiments}

In  \cite{K2}, assuming that the parameterization $\mathcal{BD}(A)$ of 
an square TP matrix $A$ is known with HRA, Plamen Koev 
presented algorithms to compute $\mathcal{BD}(A^{-1})$, the eigenvalues and
the singular values of $A$, and the solution of linear systems of equations $Ax=b$ 
where $b$ has an alternating pattern of
 signs to HRA. Koev also implemented these algorithms in order to be used
with Matlab and Octave 
in the software library {\it TNTool} available in \cite{KoevSoft}. The corresponding functions 
are \verb"TNInverse", \verb"TNEigenvalues", \verb"TNSingularValues" and
\verb"TNSolve", respectively. 
The functions require as input argument the data determining 
the bidiagonal decomposition \eqref{f} of $A$, 
 	$\mathcal{BD}(A)$ given by \eqref{eq:BD.A}, to HRA.
\verb"TNSolve" also requires a second
argument, the vector  $b$ of the linear system $Ax=b$ to be solved.
In addition, recently a function \verb"TNInverseExpand"
was added to that library, contributed by Ana Marco and Jose-Javier Mart\'inez.
This function, given $\mathcal{BD}(A)$ to HRA, returns $A^{-1}$ to HRA.

The library {\it TNTool} also provides the function 
\verb"TNProduct(B1,B2)", which, given the bidiagonal decompositions $B1$
and $B2$ to HRA of two TP matrices $F$ and $G$, provides the bidiagonal
decomposition of the TP matrix $FG$ to HRA. We can observe in the factorization 
$M=  V\,A^T$ in the proof of Theorem \ref{thm:HRA.Bes} 
that $M$ can be expressed as
the product of two TP matrices: the TP Vandermonde matrix $V$ and
the TP matrix $A^T$ defined by \eqref{eq:cambioM} and \eqref{eq:def.matC}.
Taking into account Remark \ref{rem:bd}, the bidiagonal factorization of $A^T$ to HRA
can be obtained from Theorem \ref{thm.bessel} (iii). Since $V$ is a TP Vandermonde
matrix, $\mathcal{BD}(V)$ is obtained to HRA by using \verb"TNVandBD" of library {\it TNTool}. 
Taking into account these facts,
the pseudocode providing $\mathcal{BD}(M)$ to HRA can be seen in Algorithm \ref{alg:Bes.bd}.

\begin{algorithm}[!h]
 \caption{Computation of the bidiagonal decomposition of $M$ to HRA}\label{alg:Bes.bd}
 \begin{algorithmic}
  \Require $\mathbf{t}=(t_i)_{i=0}^{n-1}$ such that $0<t_0<t_1<\ldots<t_{n-1}$
  \Ensure $B$ bidiagonal decomposition of $M$ to HRA
  \State $B1 = TNVandBD(\mathbf{t})$
  \State $B2(1,1) = 1$
  \For{$i=2:n-1$}
  	\For{$j=1:i-1$}
  		\State $B2(i,j)=\frac{(2i-2)(2i-3)}{(2i-j-1)(2i-j-2)}$
  	\EndFor
	\State $B2(i,i)=(2i-3)!!$
  	\For{$j=i+1:n-1$}
  		\State $B2(i,j)=0$
  	\EndFor
  \EndFor
  \State $B = TNProduct(B1,B2^T)$
 \end{algorithmic}
\end{algorithm}

We have implemented the previous algorithm to be used in Matlab
and Octave in a function \verb"TNBDBessel". 
The bidiagonal decompositions with HRA
obtained with this function can be used with
\verb"TNInverseExpand", \verb"TNEigenValues", \verb"TNSingularValues" and
\verb"TNSolve" in order to obtain accurate solutions for the above mentioned
algebraic problems.
Now we include some numerical tests illustrating the high 
accuracy of the new methods in contrast to the accuracy 
of the usual methods.

First we have considered the Bessel matrix of order $20$, $M_{20}$, corresponding to
the collocation matrix of the Bessel polynomials of degree at most $19$ at points
$1,2,\ldots,20$. We have computed with MATLAB the bidiagonal decomposition of $M_{20}$
to HRA with the function \verb"TNBDBessel".  Then, using this bidiagonal decomposition,
we have computed approximations to its eigenvalues and its singular values with \verb"TNEigenValues"
and \verb"TNSingularValues", respectively. We have also computed approximations to
the eigenvalues and singular values with the MATLAB functions \verb"eig" and \verb"svd",
respectively. Then we have also computed the eigenvalues and singular values of $M_{20}$
with a precision of $100$ digits using Mathematica. Taking as exact the eigenvalues and the
singular values obtained with Mathematica, we have computed the relative errors for the
approximation to the eigenvalues (resp. singular values) obtained by both \verb"eig" (resp., \verb"svd")
and \verb"TNEigenValues" (resp., \verb"TNSingularValues").

Since a Bessel matrix is STP, by Theorem 6.2 of  \cite{A} all its eigenvalues are real, distinct and positive. Taking into account this fact, the eigenvalues of $M_{20}$ have been ordered as $\lambda_1>\lambda_2>\cdots>\lambda_{20}>0$. The relative errors for the approximations to these eigenvalues of $M_{20}$ can be seen in Table \ref{tab:eig.B}. We can observe in this table that 
the approximations obtained by using the bidiagonal decomposition are very accurate. In contrast, the approximations
obtained with \verb"eig" MATLAB function are only accurate for the larger eigenvalues of $M_{20}$. In fact, the approximations
to the eigenvalues of $M_{20}$ obtained with \verb"eig" are not even positive for the smaller ones.
\begin{center}
\begin{table}[!h]
\begin{tabular}{|c|c|c|c|}
\hline
$i$ & $\lambda_i$ & rel. errors with HRA & rel. errors with \verb"eig" \\
\hline
$1$ & $4.5222e+46$ &  $4.4851e-16$ &  $1.1213e-16$ \\
$2$ & $1.2183e+42$ &  $2.5402e-16$ &  $8.0016e-15$ \\
$3$ & $7.7264e+37$ &  $2.4448e-16$ &  $3.3983e-14$ \\
$4$ & $8.7322e+33$ &  $0$                 &  $2.1003e-11$ \\
$5$ & $1.5801e+30$ &  $7.1256e-16$ &  $5.3330e-10$ \\
\vdots & \vdots &  \vdots                 &  \vdots \\
$17$ & $1.1529e+00$  & $0$  & $2.4583e+10$ \\
$18$ & $1.3072e-01$  &  $2.1233e-16$  & $2.2085e+13$ \\
$19$ & $6.1386e-03$  & $4.2389e-16$   & $1.8425e+22$ \\
$20$ & $1.2006e-04$  & $3.3864e-16$   & $4.2167e+28$ \\
 \hline
 \end{tabular}
 \caption{Relative errors for the eigenvalues of the Bessel matrix $M_{20}$}\label{tab:eig.B}
 \end{table}
 \end{center}
 
% \begin{center}
%\begin{table}[!h]
%\begin{tabular}{|c|c|c|c|}
%\hline
%$i$ & $\lambda_i$ & rel. errors with HRA & rel. errors with \verb"eig" \\
%\hline
%$1$ & $4.5222e+46$ &  $4.4851e-16$ &  $1.1213e-16$ \\
%$2$ & $1.2183e+42$ &  $2.5402e-16$ &  $8.0016e-15$ \\
%$3$ & $7.7264e+37$ &  $2.4448e-16$ &  $3.3983e-14$ \\
%$4$ & $8.7322e+33$ &  $0$                 &  $2.1003e-11$ \\
%$5$ & $1.5801e+30$ &  $7.1256e-16$ &  $5.3330e-10$ \\
%$6$ & $4.3435e+26$ &  $0$                 &  $1.3130e+01$ \\
%$7$ & $1.7656e+23$ & $3.8009e-16$ &  $2.4591e+03$ \\
%$8$ & $1.0484e+20$ & $1.5628e-16$ & $1.6828e+03$ \\
%$9$ & $9.0766e+16$ & $1.7628e-16$ & $1.1853e+03$ \\
%$10$ & $1.1534e+14$ &          $0$  & $1.4105e+03$ \\
%$11$ & $2.1833e+11$ & $6.9888e-16$  & $7.0247e+04$ \\
%$12$ & $6.3074e+08$  & $3.7800e-16$  & $1.9551e+06$ \\
%$13$ & $2.8844e+06$  & $4.8432e-16$  & $7.6261e+05$ \\
%$14$ & $2.2083e+04$  & $8.2369e-16$  & $7.0232e+04$ \\
%$15$ & $3.1128e+02$  & $9.1307e-16$ & $2.3765e+06$ \\
%$16$ & $1.0055e+01$   & $3.5332e-16$ & $9.0055e-01$ \\
%$17$ & $1.1529e+00$  & $0$  & $2.4583e+10$ \\
%$18$ & $1.3072e-01$  &  $2.1233e-16$  & $2.2085e+13$ \\
%$19$ & $6.1386e-03$  & $4.2389e-16$   & $1.8425e+22$ \\
%$20$ & $1.2006e-04$  & $3.3864e-16$   & $4.2167e+28$ \\
% \hline
% \end{tabular}
% \caption{Relative errors for the eigenvalues of the Bessel matrix $M_{20}$}\label{tab:eig.B}
% \end{table}
% \end{center}
 
The $20$ real and positive singular values of $M_{20}$ have also been ordered as $\sigma_1\ge\sigma_2\ge\cdots\ge\sigma_{20}>0$. The relative errors for the approximations to these singular values of $M_{20}$ can be seen in Table \ref{tab:svd.B}.
As in the case of the eigenvalues, the approximation obtained for all the singular values using the bidiagonal decomposition 
are very accurate, but only the approximations obtained for the larger singular values by using \verb"svd" are accurate.
 
 \begin{center}
\begin{table}[!h]
\begin{tabular}{|c|c|c|c|}
\hline
$i$ & $\sigma_i$ & rel. errors with HRA & rel. errors with \verb"svd" \\
\hline
$1$  & $4.8763e+46$ &  $2.0797e-15$  & $4.1594e-16$ \\
$2$  & $1.5204e+42$ &  $6.1065e-16$ &  $4.0710e-16$ \\
$3$  & $1.1076e+38$ &  $3.4108e-16$  & $9.1922e-14$ \\
$4$  & $1.4266e+34$ &  $8.0818e-16$  & $8.3458e-11$ \\
$5$   & $2.9165e+30$ & $3.8604e-16$ & $1.2661e-03$ \\
\vdots & \vdots & \vdots &  \vdots \\
$17$ & $1.0795e+00$ &  $4.1139e-16$ &  $5.0876e+11$ \\
$18$ & $1.5106e-02$   & $2.1818e-15$ &  $3.4860e+12$ \\
$19$ & $9.1285e-05$ &  $2.0785e-15$ &  $1.2718e+13$ \\
 $20$ & $1.6258e-07$ & $3.2563e-16$ &  $1.9230e+06$ \\
    \hline
 \end{tabular}
 \caption{Relative errors for the singular values of the Bessel matrix $M_{20}$}\label{tab:svd.B}
 \end{table}
 \end{center}
 
%\begin{center}
%\begin{table}[!h]
%\begin{tabular}{|c|c|c|c|}
%\hline
%$i$ & $\sigma_i$ & rel. errors with HRA & rel. errors with \verb"eig" \\
%\hline
%$1$  & $4.8763e+46$ &  $2.0797e-15$  & $4.1594e-16$ \\
%$2$  & $1.5204e+42$ &  $6.1065e-16$ &  $4.0710e-16$ \\
%$3$  & $1.1076e+38$ &  $3.4108e-16$  & $9.1922e-14$ \\
%$4$  & $1.4266e+34$ &  $8.0818e-16$  & $8.3458e-11$ \\
%$5$   & $2.9165e+30$ & $3.8604e-16$ & $1.2661e-03$ \\
%$6$  & $8.9750e+26$ &  $7.6568e-16$ &  $9.3463e+02$ \\
%$7$  & $4.0432e+23$ &  $1.6598e-16$ &  $2.1085e+03$ \\
%$8$   & $2.6321e+20$ &  $1.2450e-15$ & $3.6454e+05$ \\
%$9$  & $2.4690e+17$ &           $0$ &  $1.3831e+06$ \\
%$10$ & $3.3556e+14$ &  $3.7251e-16$  & $2.6205e+07$ \\
%$11$ & $6.6978e+11$ & $5.4676e-16$ &  $3.4324e+08$ \\
%$12$ &  $2.0077e+09$ &  $3.5626e-16$ & $1.8733e+08$ \\
%$13$ & $9.3469e+06$ &  $1.1957e-15$ &  $1.9697e+10$ \\
%$14$ & $7.1062e+04$ &  $4.0955e-16$ &  $4.0300e+09$ \\
%$15$ & $9.4908e+02$ &  $3.5936e-16$ &  $9.2542e+10$ \\
%$16$ & $2.7180e+01$ &  $6.5356e-16$ &  $6.5490e+11$ \\
%$17$ & $1.0795e+00$ &  $4.1139e-16$ &  $5.0876e+11$ \\
%$18$ & $1.5106e-02$   & $2.1818e-15$ &  $3.4860e+12$ \\
%$19$ & $9.1285e-05$ &  $2.0785e-15$ &  $1.2718e+13$ \\
% $20$ & $1.6258e-07$ & $3.2563e-16$ &  $1.9230e+06$ \\
%    \hline
% \end{tabular}
% \caption{Relative errors for the singular values of the Bessel matrix $M_{20}$}\label{tab:svd.B}
% \end{table}
% \end{center}
 
 We have also computed approximations to the inverse $(M_{20})^{-1}$ with both \verb"inv" and
 \verb"TNInverseExpand" using the bidiagonal decomposition provided by \verb"TNBDBessel".
 We have also obtained with Mathematica the inverse using exact arithmetic. Then,
 we have computed
 the componentwise relative errors for both approximations. In Table \ref{tab:inv.B}
 the mean and maximum componentwise relative errors are shown. It can be seen that
 the inverse obtained with the HRA methods is very accurate in contrast to the poor 
 approximation obtained with \verb"inv".

\begin{center}
\begin{table}[!h]
\begin{tabular}{|r|c|c|}
\hline
 &  rel. errors with HRA & rel. errors with \verb"inv" \\
\hline
mean &   $1.8498e-16$ & $2.1364e-01$ \\
maximum & $8.4304e-16$ & $3.0878e-01$ \\
      \hline
 \end{tabular}
 \caption{Relative errors for the inverse  of the Bessel matrix $M_{20}$}\label{tab:inv.B}
 \end{table}
 \end{center}
 
 We have considered two systems of linear equations
 \[
 	M_{20}x=b^1\quad\text{and}\quad M_{20}x=b^2,
 \]
 where the entries of $b^2$ are randomly generated as integers in the interval $[1,1000]$ and the $i$-th
 entry of $b^1$ is given by $b^1_{i}=(-1)^{i+1}b^2_i$ for $i=1,\ldots,20$. So, the independent vector of 
 the system $M_{20}x=b^1$ has an alternating pattern of signs and the linear system can be solved with HRA by Theorem \ref{thm:HRA.Bes}. Approximations $\widehat{x}$ to the 
 solutions $x$ of both linear systems have been obtained with MATLAB, the first
one using \verb"TNSolve" and the bidiagonal decomposition of the Bessel matrix obtained with \verb"TNBDBessel", and
the second one using the usual MATLAB command \verb"A\b". By using Mathematica with exact arithmetic, the
exact solution of the systems have been computed, and then, the componentwise relative errors for the two approximations
obtained with MATLAB have been computed. Table \ref{eq:sis1.BV} shows the componentwise relative errors corresponding to the 
system $M_{20}x=b^1$. It can be observed that the approximation to the solution provided by \verb"TNSolve" is
very accurate. This fact can be expected since the independent vector $b^1$ of the system has an alternating pattern of signs
and then it is known that \verb"TNSolve" provides the solution to HRA (see \cite{K2}). On the other hand, it can
also be observed in the table that the accuracy
of the approximation provided by \verb"A\b" is very poor.

 \begin{center}
\begin{table}[!h]
\begin{tabular}{|r|c|c|}
\hline
$i$ &  $\left|\frac{\widehat{x_i}-x_i}{x_i}\right|$ with HRA & $\left|\frac{\widehat{x_i}-x_i}{x_i}\right|$ with \verb"A\b" \\
\hline
 $1$ & $5.6243e-16$ & $1.8543e-01$ \\
 $2$ & $1.8069e-16$ & $1.8643e-01$ \\
 $3$ & $0$ & $1.8822e-01$ \\
 $4$ & $1.2831e-16$ & $1.9058e-01$ \\
 $5$ & $2.2120e-16$ &  $1.9336e-01$ \\
 \vdots & \vdots &  \vdots \\
$17$ & $0$ & $2.3342e-01$ \\
$18$ & $1.4910e-16$ & $2.3627e-01$ \\
$19$ & $2.3844e-16$ & $2.3897e-01$ \\
$20$ &  $0$ &    $2.4153e-01$ \\
\hline
 \end{tabular}
 \caption{Relative errors for the solution of the linear system $M_{20}x=b^1$}\label{eq:sis1.BV}
 \end{table}
 \end{center}

%  \begin{center}
%\begin{table}[!h]
%\begin{tabular}{|r|c|c|}
%\hline
%$i$ &  $\left|\frac{\widehat{x_i}-x_i}{x_i}\right|$ with HRA & $\left|\frac{\widehat{x_i}-x_i}{x_i}\right|$ with \verb"A\b" \\
%\hline
% $1$ & $5.6243e-16$ & $1.8543e-01$ \\
% $2$ & $1.8069e-16$ & $1.8643e-01$ \\
% $3$ & $0$ & $1.8822e-01$ \\
% $4$ & $1.2831e-16$ & $1.9058e-01$ \\
% $5$ & $2.2120e-16$ &  $1.9336e-01$ \\
% $6$ & $3.2671e-16$ &  $1.9644e-01$ \\
% $7$ & $0$ & $1.9973e-01$ \\
% $8$ & $1.7213e-16$ & $2.0317e-01$ \\
% $9$ & $2.2534e-16$ & $2.0669e-01$ \\
%$10$ & $0$ & $2.1025e-01$ \\
%$11$ & $0$ &  $2.1380e-01$ \\
%$12$ & $1.2441e-16$ & $2.1730e-01$ \\
%$13$ & $1.5484e-16$ & $2.2074e-01$ \\
%$14$ & $2.6516e-16$ & $2.2409e-01$ \\
%$15$ & $0$ & $2.2732e-01$ \\
%$16$ & $1.3155e-16$ & $2.3044e-01$ \\
%$17$ & $0$ & $2.3342e-01$ \\
%$18$ & $1.4910e-16$ & $2.3627e-01$ \\
%$19$ & $2.3844e-16$ & $2.3897e-01$ \\
%$20$ &  $0$ &    $2.4153e-01$ \\
%\hline
% \end{tabular}
% \caption{Relative errors for the solution of the linear system $M_{20}x=b_1$}\label{eq:sis1.BV}
% \end{table}
% \end{center}

Table \ref{eq:sis2.BV} shows the componentwise relative errors corresponding to the 
system $M_{20}x=b^2$. In this case, since the independent vector $b^2$ has not an 
alternating pattern of signs, it is not guaranteed to obtain an approximation to HRA
by using \verb"TNSolve". However, it can be observed that in this case the approximation to the 
solution provided by \verb"TNSolve" is also very accurate in contrast to the poor accuracy of
the approximation provided by \verb"A\b". 
\begin{center}
\begin{table}[!h]
\begin{tabular}{|r|c|c|}
\hline
$i$ &  $\left|\frac{\widehat{x_i}-x_i}{x_i}\right|$ with \verb"TNSolve" & $\left|\frac{\widehat{x_i}-x_i}{x_i}\right|$ with \verb"A\b" \\
\hline
$1$ & $1.6814e-16$  &  $8.9852e-01$ \\
$2$ &   $1.2521e-16$  &  $8.3698e-01$ \\ 
$3$ &  $1.7413e-16$  & $7.4661e-01$ \\
$4$ &   $2.7288e-16$ &  $6.5449e-01$ \\
$5$ &    $2.0323e-16$ &  $5.7258e-01$ \\
\vdots &   \vdots &  \vdots \\
$17$ &   $1.8808e-16$  & $2.1570e-01$ \\
$18$ &   $1.6663e-16$  & $2.0699e-01$ \\
$19$ &   $2.5411e-16$  & $1.9933e-01$ \\
$20$ &        $0$  & $1.9257e-01$ \\
\hline
 \end{tabular}
 \caption{Relative errors for the solution of the linear system $M_{20}x=b^2$}\label{eq:sis2.BV}
 \end{table}
 \end{center}

% \begin{center}
%\begin{table}[!h]
%\begin{tabular}{|r|c|c|}
%\hline
%$i$ &  $\left|\frac{\widehat{x_i}-x_i}{x_i}\right|$ with \verb"TNSolve" & $\left|\frac{\widehat{x_i}-x_i}{x_i}\right|$ with \verb"A\b" \\
%\hline
%$1$ & 1.6814e-16  &  8.9852e-01 \\
%$2$ &    1.2521e-16  &  8.3698e-01 \\ 
%$3$ &   1.7413e-16  & 7.4661e-01 \\
%$4$ &    2.7288e-16 &  6.5449e-01 \\
%$5$ &    2.0323e-16 &  5.7258e-01 \\
%$6$ &    1.3026e-16 &  5.0374e-01 \\
%$7$ &   2.6862e-16 &  4.4715e-01 \\
%$8$ &            0  & 4.0090e-01 \\
%$9$ &            0  &  3.6303e-01 \\
%$10$ &            0 &  3.3183e-01 \\
%$11$ &            0  & 3.0593e-01 \\
%$12$ &   2.0526e-16 &  2.8427e-01 \\
%$13$ &   2.3582e-16 &  2.6599e-01 \\
%$14$ &            0 &  2.5046e-01 \\
%$15$ &            0  & 2.3715e-01 \\
%$16$ &            0  & 2.2567e-01 \\
%$17$ &   1.8808e-16  & 2.1570e-01 \\
%$18$ &   1.6663e-16  & 2.0699e-01 \\
%$19$ &   2.5411e-16  & 1.9933e-01 \\
%$20$ &             0  & 1.9257e-01 \\
%\hline
% \end{tabular}
% \caption{Relative errors for the solution of the linear system $M_{20}x=b_2$}\label{eq:sis2.BV}
% \end{table}
% \end{center}             

It can be observed that the smaller an eigenvalue (resp., singular value) is,
the larger the relative error corresponding to the usual methods is. So, now let us consider 
the Bessel matrices $M_n$ of order $n$, for $n=2,\ldots,15 $ given by the 
collocation matrices
of the Bessel polynomials $(B_0(x),\ldots,B_{n-1}(x))$ at 
the points $1,\ldots,n$, that is, $M_n=(B_{j-1}(i))_{1\leq i,j\leq n}$.
In the same way that in the previous examples we have computed the eigenvalues,
the singular values and the inverses of these matrices both with the usual MATLAB
functions and to HRA by using \verb"TNBDBessel". Then we have computed the relative
errors for the approximation to the smallest eigenvalue and the smallest singular value of 
each matrix, and the componentwise relative error for the approximations to the inverses.

\begin{figure}
  \includegraphics[scale=0.3]{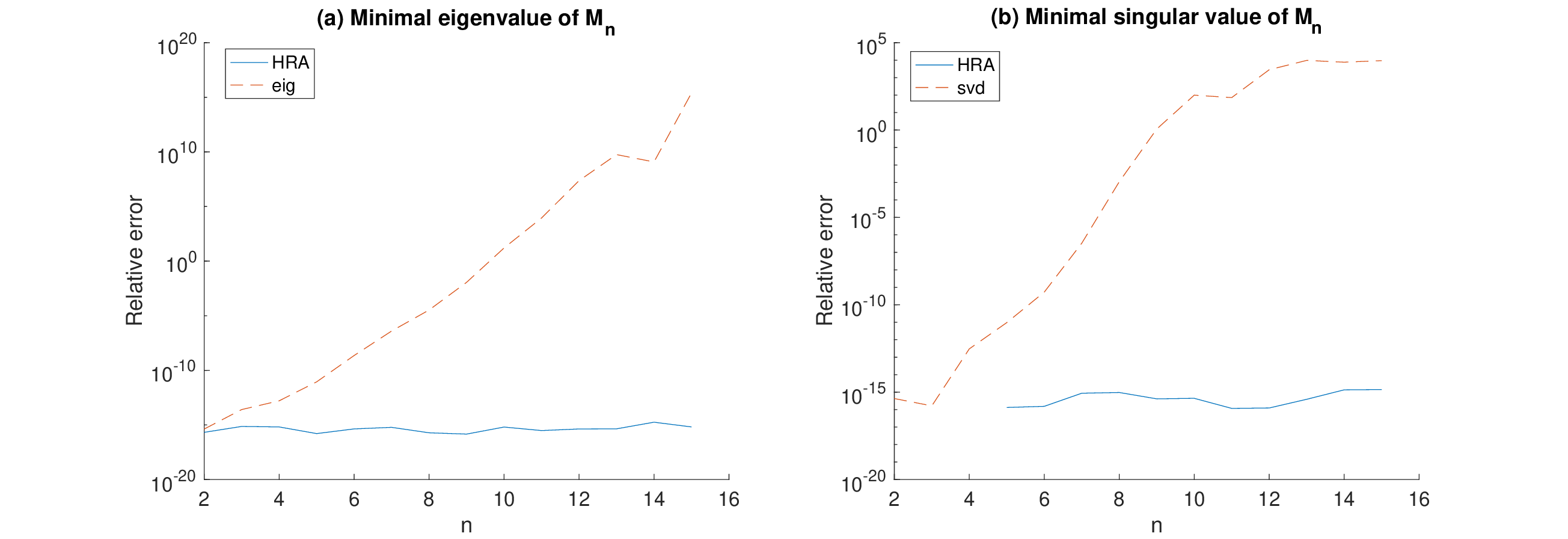}
\caption{Relative error for the minimal eigenvalue and singular value of $M_n$}
\label{fig:val}       % Give a unique label
\end{figure}

The relative errors for the smallest eigenvalues and the smallest singular values of the
Bessel matrices $M_n$, $n=2,\ldots,15$, can be seen in Figure \ref{fig:val} (a) and 
(b), respectively.

\begin{figure}
  \includegraphics[scale=0.3]{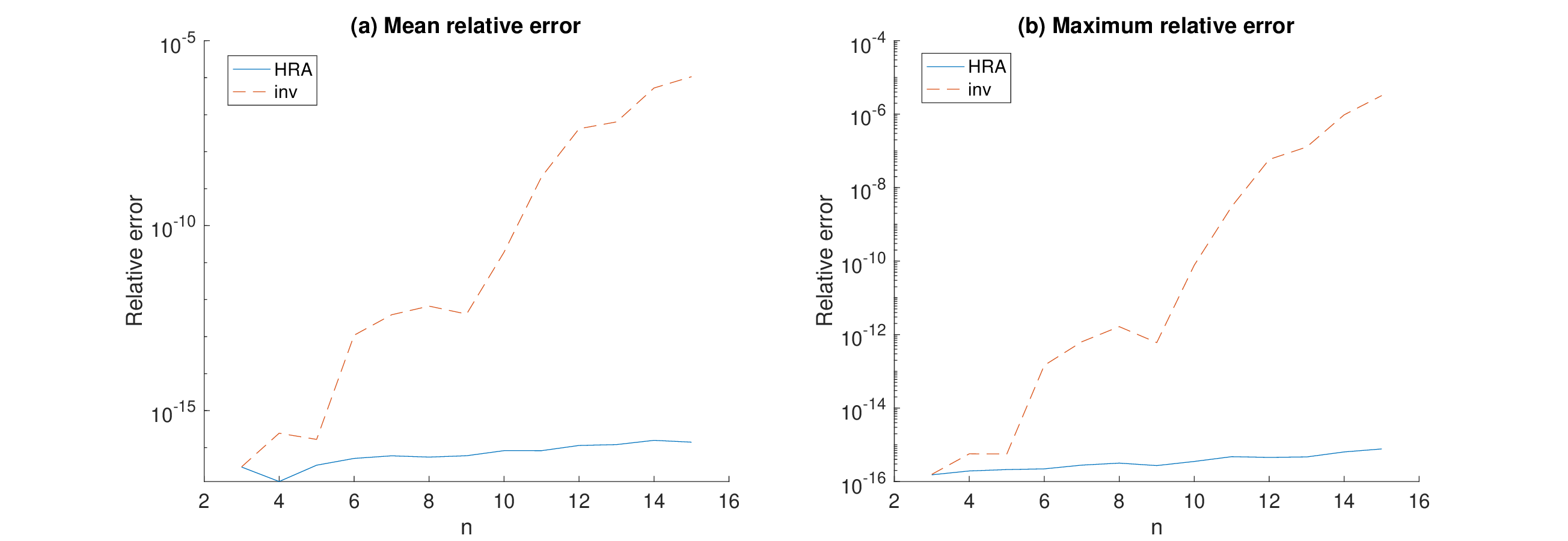}
\caption{Relative errors for $(M_n)^{-1}$}
\label{fig:inv}       % Give a unique label
\end{figure}

The mean and the maximum componentwise relative errors corresponding 
to the approximation of the inverses $(M_n)^{-1}$ can be seen in Figure
\ref{fig:inv} (a) and (b), respectively.

In an analogous way to the Bessel matrix we can derive an algorithm to obtain
the bidiagonal decomposition of a reverse Bessel matrix to HRA. So,
the pseudocode providing $\mathcal{BD}(M_r)$ to HRA can be seen in Algorithm \ref{alg:rev.Bes.bd}.

\begin{algorithm}[!h]
 \caption{Computation of the bidiagonal decomposition of $M_r$ to HRA}\label{alg:rev.Bes.bd}
 \begin{algorithmic}
  \Require $\mathbf{t}=(t_i)_{i=0}^{n-1}$ such that $0<t_0<t_1<\ldots<t_{n-1}$
  \Ensure $B$ bidiagonal decomposition of $M_r$ to HRA
  \State $B1 = TNVandBD(\mathbf{t})$
  \For{$i=1:n-1$}
  	\For{$j=1:i-1$}
		\If{$j$ is odd}
  			\State $B2(i,j)=2i-2j-1$
		\EndIf
  	\EndFor
	\State $B2(i,i)=1$
  	\For{$j=i+1:n-1$}
  		\State $B2(i,j)=0$
  	\EndFor
  \EndFor
  \State $B = TNProduct(B1,B2^T)$
 \end{algorithmic}
\end{algorithm}

\begin{figure}
  \includegraphics[scale=0.3]{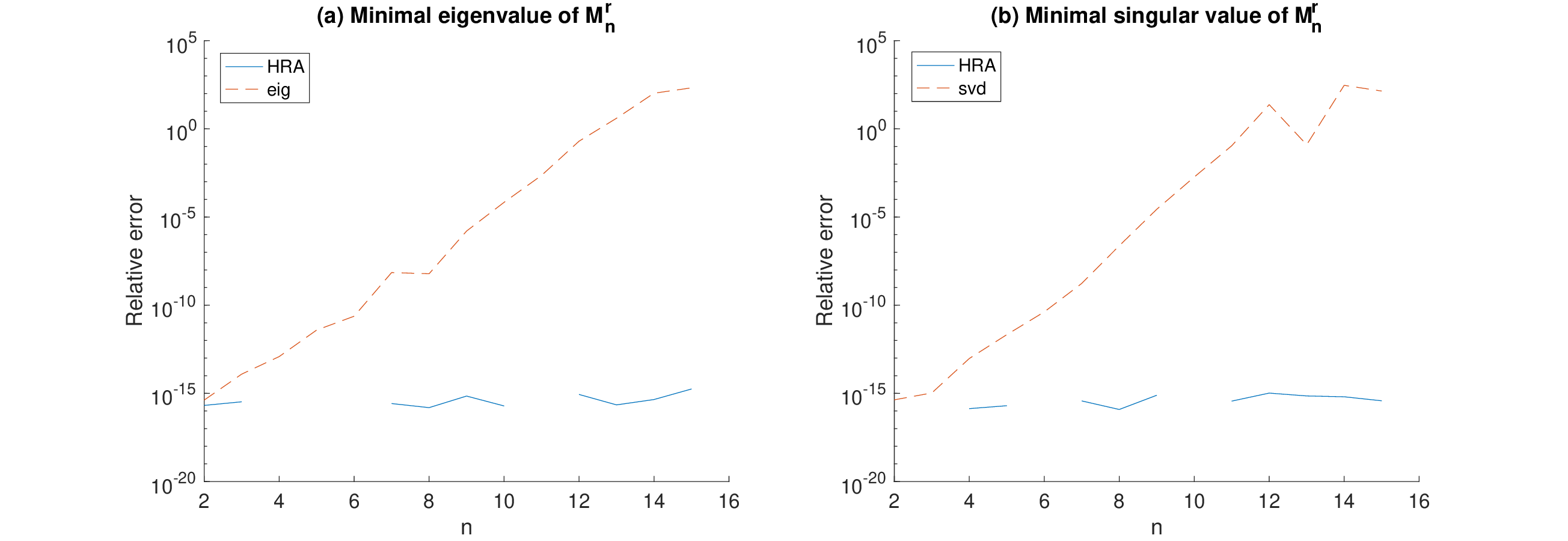}
\caption{Relative error for the minimal eigenvalue and singular value of $M_n^r$}
\label{fig:valR}       % Give a unique label
\end{figure}

\begin{figure}
  \includegraphics[scale=0.3]{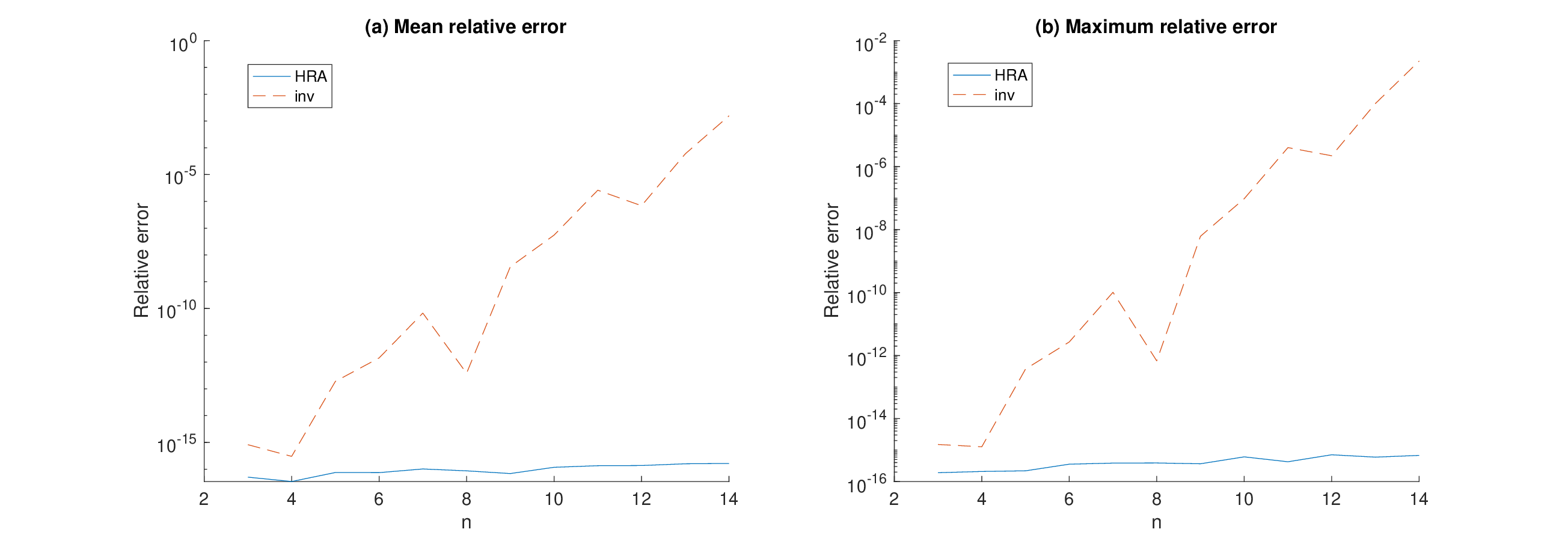}
\caption{Relative errors for $(M_n^r)^{-1}$}
\label{fig:invR}       % Give a unique label
\end{figure}

For the reverse Bessel matrices we have carried out the same numerical tests as for the Bessel matrices
and we have deduced exactly the same conclusions. For the sake of brevity, for the reverse Bessel matrices
only the relative errors for the smallest eigenvalue and singular value, and the componentwise mean and 
maximum relative error for the inverses of the reverse Bessel matrices $M_n^r=(B_{j-1}^r(i))_{1\leq i,j\leq n}$,
$n=2,\ldots,15$, are shown in Figures \ref{fig:valR} and \ref{fig:invR}, respectively.

\begin{acknowledgements}
This work was partially supported through the Spanish research grant MTM2015-65433-P (MINECO/FEDER), by
Gobierno de Aragón, and by Fondo Social Europeo.
\end{acknowledgements}

\end{document}